\input{style/arxiv-ba.cfg}
\documentclass[ba,linksfromyear,preprint]{imsart}
\usepackage{amssymb}
\usepackage{vtexurl}

\pubyear{2015}
\volume{10}
\issue{1}
\firstpage{227}
\lastpage{231}
\doi{10.1214/14-BA936}% Updated by VTEXPTS2LaTeX.exe, 15.01.2015 15:01

\startlocaldefs
\newcommand{\enquote}[1]{``#1''}
\newcommand{\bm}[1]{\mbox{\boldmath$#1$}}

\let\bid@start\@empty
\let\bid@end\@empty

\def\MR@url{http://www.ams.org/mathscinet-getitem?mr=}
\def\MR#1{\href{\MR@url#1}{MR#1}}
\def\BDOI#1{%
\edef\doi@base@i{\doi@base}\def\doi@base{}%
doi:~\doiurl{\doi@base@i#1}}

\appto\bid@start{\def\doi@size{\ttfamily}}
\appto\bid@end{\unskip.}

\define@key{bid}{mr}{\def\bid@mr{\MR{#1}}}
\define@key{bid}{doi}{\def\bid@doi{\BDOI{#1}}}
\define@key{bid}{pubmed}{}
\define@key{bid}{pii}{}
\define@key{bid}{issn}{}

\def\bid#1{%
       \bgroup
       \bid@start
       \let\bid@output\@empty
       \setkeys{bid}{#1}\ignorespaces%
       \ifdefvoid\bid@mr{}{\appto\bid@output{\bid@mr}}%
       \ifdefvoid\bid@doi{}{
         \ifdefempty\bid@output{}{\appto\bid@output{. }}% separator if output not empty
         \appto\bid@output{\bid@doi}%
       }%
       \bid@output
       \bid@end
       \egroup
}
\endlocaldefs

\begin{document}

\begin{frontmatter}
\title{Comment on Article by Berger, Bernardo, and~Sun\thanksref{T1}}
\runtitle{Comment on Article by Berger, Bernardo, and Sun}

\relateddois{T1}{Main article DOI: \relateddoi[ms=BA915,title={James O. Berger, Jose M. Bernardo, and Dongchu Sun. Overall Objective Priors}]{Related item:}{10.1214/14-BA915}.}

\begin{aug}
\author[a]{\fnms{Manuel} \snm{Mendoza}\ead[label=e1]{mendoza@itam.mx}}
\and
\author[b]{\fnms{Eduardo} \snm{Guti\'errez-Pe\~na}\corref{}\ead[label=e2]{eduardo@sigma.iimas.unam.mx}}

\runauthor{M. Mendoza and E. Guti\'errez-Pe\~na}

\address[a]{Dept.\ of Statistics, ITAM, Mexico, \printead{e1}}
\address[b]{Dept.\ of Probability and Statistics, UNAM, Mexico,
\printead{e2}}
\end{aug}

\end{frontmatter}

%% Mainmatter %%

%s1 ###
\section{Introduction}

The search of a prior distribution $p(\bm{\omega})$ to be used as part
of an objective Bayesian analysis of a model $p(x | \bm{\omega})$ has
proved to be a formidable endeavour. This is an area where we do not
have a definitive answer yet, and any contribution to the understanding
of the subject must be welcome. The authors of this paper are among the
most prominent contributors to this field, and reading the manuscript
has been very stimulating.

Research on the problem has mainly dealt with three issues: first, a
definition of what a {\em non-informative, reference} or {\em
objective} prior $p(\bm{\omega})$ must be; second, an operational
algorithm to calculate such priors; third, the evaluation of the
resulting prior(s) in accordance to certain criteria such as
invariance, the avoidance of paradoxes, or desirable frequentist properties.

To us, and this is a subjective judgment, the most convincing approach to
produce this sort of priors is reference analysis
\citep{Bernardo1979,BerBer1992a,BerBer1992b,Bernardo2005,Bergeretal2009}. This
procedure: (i) defines the reference prior as the prior maximizing the
expected gain of information provided by a sample; (ii) includes a
general (although potentially involved) algorithm to calculate the
prior; and (iii)~avoids a number of paradoxes. Moreover, it generalizes
the Jeffreys prior and exhibits its limitations. Among its most
remarkable results, it shows that the form of the reference prior $p(\bm
{\omega})$ may depend on the function of the parameters $\bm{\theta} =
\bm{\theta}(\bm{\omega})$ which is considered by the researcher to be
of main interest.

Since its inception, the algorithm to obtain reference priors has
evolved. This is the case specifically in the multiparameter setting.
The most recent version \citep{BerBer1992a,BerBer1992b} requires all
scalar components of the parameter to be {\em strictly} ordered in
terms of their inferential interest. Thus, in principle, the current
approach does not offer any solution if the researcher is
simultaneously interested in two or more scalar parameters (or
functions thereof). Interestingly, the original algorithm of \cite{Bernardo1979} did cover this situation, although the solution was the
multivariate Jeffreys prior which leads to unsettling paradoxes in some cases.

In this paper, the authors explore some ideas to extend the reference
analysis to this yet unsolved case. They also seem to be considering a
more general version of the problem by assuming that the number of
scalar parameters (or functions of the parameters) of interest may be
greater than the number of parameters in the model. So, the question
is: What should the objective prior $\pi^R(\bm{\omega})$ ($\bm{\omega}
\in\mathbb{R}^k$) be if there are $m$ functions ($\theta_1(\bm{\omega
}), \theta_2(\bm{\omega}),\dots,\theta_m(\bm{\omega})$) which are of
simultaneous interest, where $m$ is not constrained to be less than or
equal to $k$? Three methods to produce the required prior distribution
are discussed: (i) the common reference prior; (ii) the reference
distance approach; and (iii) the hierarchical approach.

%s2 ###
\section{Common reference prior}

This is not really a method. If the reference priors corresponding to
$\theta_i(\bm{\omega})$ as the parameter of interest ($i=1,\ldots, m$)
are the same for any ordering of the remaining parameters, then the
posed problem simply vanishes. It is interesting to see some examples
illustrating particular cases where the common prior exists, but it is
desirable -- and would be much more useful -- to have general results
characterizing sampling models where, for example, Theorem 2.1 applies
and hence a common reference prior may be found. In this regard,
results such as those in \citet{gpr} and \citet{cvgp} could provide a
good starting point. These authors find reference priors for wide
classes of exponential families that include the family discussed in
Section 2.1.3 of the present paper as a particular case.

It must be pointed out that this section relies on the analysis of the
information matrix $\bm{I}(\bm{\theta})$, so all reviewed scenarios
assume $m \leq k$. Also, a somewhat disquieting result is that of
Section 2.2.2, where the authors show that
$\pi^R(\psi_1,\psi_2,\psi_3,\mu_1,\mu_2) \propto(\psi_1\psi_2)^{-1}$
is the one-at-a-time reference prior for any of these parameters and
any possible ordering. In particular, it is the reference prior for the
case where $\mu_2$ is the parameter of main interest. It so happens,
however, that this prior is equivalent to the right-Haar prior which
leads to a problematic posterior precisely for $\mu_2$. This result
would imply that, in general, reference analysis might produce
inadequate posteriors {\em for the parameter of interest}, depending on
the specific accompanying parameters.

%s3 ###
\section{Reference distance method}

In order to introduce this method, the authors explicitly assume that
$\bm{\theta} = \bm{\omega}$, hence $m = k$. The idea is to find an
overall prior $\pi(\bm{\theta})$ such that each of its marginal
posteriors $\pi(\theta_i|\bm{x})$ is close to the corresponding
marginal posterior $\pi_i(\theta_i|\bm{x})$ obtained when $\theta_i$ is
the parameter of interest ($i=1,\ldots, m$). As a measure of
approximation the authors propose a weighted average of expected
logarithmic divergences, although other measures could in principle be
used. Also, the search for the overall prior is restricted to a
specific parametric family
$\mathcal{F} = \{\pi(\bm{\theta}| \bm{a}), \bm{a} \in\mathcal{A} \}$.
Apart from the fact (acknowledged by the authors) that the existence of
an optimal $\bm{a}$ is not guaranteed, a rather unappealing feature of
this proposal is its dependence on the family $\mathcal{F}$. The authors
offer no guidance on how to choose $\mathcal{F}$ \emph{in general}. If the
aim is to produce an objective approach, it seems desirable that $\mathcal
{F}$ be somehow \emph{intrinsic} to the sampling model. The examples in
the paper suggest that perhaps this could be achieved through some kind
of conjugacy.

Incidentally, the reference distance method bears some resemblance to
the mean-field approach to variational inference, which is relatively
straightforward in the case of exponential families with conjugate
priors; see, for example, \citet[Chapter~10]{prml}. What is the
authors' take on this?

We would like now to comment on Example 3.2.4. There, the normal model
$N(x|\mu,\sigma)$ is considered, and the parameters of interest are $\mu
$, $\sigma$ and $\phi= \mu/\sigma$. (Note that, despite the authors'
remark at the beginning of Section 3, here $\bm{\theta} \neq\bm{\omega
}$ and $m>k$.) In any case, the authors remind us that the reference
prior when $\mu$ or $\sigma$ is the parameter of interest is $\pi(\mu
,\sigma)= \sigma^{-1}$, whereas the reference prior for $\phi= \mu
/\sigma$ is given by
$\pi_{\phi}(\mu,\sigma) = (2\sigma^2+\mu^2)^{-1/2}\sigma^{-1}$. They
then propose, as a ``natural'' choice, the class of relatively
invariant priors
$\mathcal{F}= \{ \pi( \mu,\sigma)=\sigma^{-a}; a>0 \}$.
For this family, they show that the overall prior for $(\mu,\sigma, \phi
)$ can be approximated by $\pi^o(\mu,\sigma) = \sigma^{-1}$, so that
inclusion of $\phi$ as an additional parameter of interest makes no
difference. We find this rather disappointing. From an algorithmic
point of view, this outcome is not surprising given the choice of $\mathcal
{F}$ and the form of the reference priors for $\mu$ and $\sigma$. Only
a large weight on the divergence corresponding to $\phi$ could lead to
a different result. An idea that springs to mind is to try another
(arguably more ``natural'' family) such as $\mathcal{G} = \{ \pi( \mu
,\sigma| a_1,a_2)=(2\sigma^2+\mu^2)^{-a_1}\sigma^{-a_2}; a_1>0, a_2>0\}
$, which includes all three reference priors for $\mu$, $\sigma$ and~$\phi$. On the other hand, since $\pi_{\mu}(\mu,\sigma)$ and $\pi
_{\sigma}(\mu,\sigma)$ are equal in this case, the authors could
alternatively have minimized the sum of the two divergences
corresponding to the marginal posterior of $\phi$ and the {\em joint}
posterior of $(\mu,\sigma)$). We wonder how these alternative ideas
compare with that proposed in the paper for this example.

%s4 ###
\section{Hierarchical approach}

The idea of this approach is, first, to find a ``natural'' parametric
family of proper priors $\pi(\bm{\theta}|a)$ such that $a \in\mathbb
{R}$ and the integrated likelihood results in a proper density
$p(x|a)$. Then, the univariate reference prior for $a$, $\pi^R(a)$, is
obtained for this latter model. Finally, the overall prior $\pi^o(\bm
{\theta})$ is defined as the expectation of $\pi(\bm{\theta}|a)$ with
respect to $\pi^R(a)$. This is an intuitive and seemingly reasonable
idea. However, it is not clear how to make explicit that $\bm{\theta}$
is the parameter of interest even though the model is originally
indexed by $\bm{\omega}$, especially when the dimension of $\bm{\theta
}$ is larger than that of~$\bm{\omega}$. (See the comment below
concerning the multi-normal means example.) We wonder if the authors
can provide some advice on how this could be achieved in general. On
the other hand, as in the reference distance case, dependence upon a
specific family of priors introduces no small amount of arbitrariness
in the method. Here, again, a proper objective method would use an \emph
{intrinsic} family entirely determined by the sampling model. One
possibility, particularly suitable for the case of hierarchical models,
would be to elaborate on the idea of conjugate likelihood distributions
\citep{gms}, although a suitable restriction should be imposed on the
corresponding conjugate family in order to get a one-dimensional
hyperparameter. Concerning the implementation of the method, the
authors suggest that integration to get the overall prior can be
avoided by using $\pi^o(\bm{\theta}) = \pi(\bm{\theta}|\hat{a})$
instead, where $\hat{a}$ is the mode of the posterior $p(a|\bm{x})$.
This proposal may be efficient from a computational point of view, but
it is both surprising and disappointing since it essentially reduces
the hierarchical approach to a standard empirical Bayes procedure and
leads to a data-dependent prior.

The example in Section 4.2 concerning the multivariate hypergeometric
model is confusing and does not quite illustrate the method described
above. First, the parameters of the sampling model are given a
multinomial prior (which does not depend on a single scalar parameter
$a$, but on a vector of probabilities $\bm{p}_k$); then, the likelihood
is integrated and shown to yield a multinomial distribution. In the
process, the $k$ original parameters $R_1,R_2,\ldots,R_k$ are replaced
by the parameters $p_1,p_2,\ldots,p_k$, so the idea of reducing the
problem to the determination of the reference prior for a scalar
parameter is abandoned. Next, in the multinomial model, the approximate
overall prior obtained {\em using the reference distance method} is
adopted for the hyperparameters $\bm{p}_k$. Finally, the corresponding
integrated Multinomial--Dirichlet distribution is declared as the
overall prior for $R_1,R_2,\ldots,R_k$. We find this \emph{ad hoc}
combination of methods difficult to justify as a general procedure.

An alternative formulation could be based on the idea of
super-populations (quite common in the field of survey sampling) as
follows. Let us assume that a random sample of size $N$ is obtained
from a multinomial distribution $Mu_k(\bm{Y}_k|1,\bm{p}_k)$. As a
result we get a vector $R_1,R_2,\ldots,R_k$ describing the number of
sampled units in each category. Now imagine that we then get a
subsample of size $n$, {\em without replacement}, from the sample of
size $N$. In this setting, the multinomial distribution describes an
infinite super-population, the sample of size $N$ is the finite
population of interest and the subsample of size $n$ is the actual
sample we observe. Given the sample, the likelihood based on the
subsample corresponds to that of a hypergeometric distribution.
However, with respect to the super-population, the subsample is just a
sample of the original multinomial population whose parameters are
given by the vector $\bm{p}_k$. Within this framework, $R_1,R_2,\ldots
,R_k$ are observables and any inference regarding these quantities must
be produced through the corresponding posterior predictive
distribution. This argument shows that the hypergeometric problem can
be viewed as a multinomial one where the interest is not really on the
parameters but on observables, and the relevant overall prior is that
for $\bm{p}_k$, no matter which method we use.

The example on the multi-normal means (Section 4.3) deserves a few
words as well. Here, the parameters of interest are, using the same
notation as the authors, $\mu_i$; $i=1,\ldots,m$ and $|\bm{\mu}|^2 =
\mu_1^2 +\cdots+\mu_m^2$. First, we note that throughout the paper $k$
refers to the dimension of $\bm{\omega}$ and $m$ is the number of
parameters of interest (the dimension of~$\bm{\theta}$), so in this
example we have $m=k+1$. It must be pointed out, however, that the
hierarchical method, as defined, cannot be applied when $m>k$ since the
distribution $\pi(\bm{\theta}|a)$ would then be defined over a space of
functionally related components of $\bm{\theta}$ and would be singular.
This fact is implicitly recognized by the authors when they propose a
prior for $(\mu_1,\ldots,\mu_m)$ only, ignoring the last parameter of
interest, $|\bm{\mu}|^2$. They then argue that the resulting overall
prior is reasonable not only for each mean $\mu_i$ but also for $|\bm
{\mu}|^2$. The key issue here is the convenient choice of $\pi(\bm{\mu
}|a)$ as the product of the normals $N(\mu_i|0,a)$. So, strictly
speaking, this problem is not actually solved by using the hierarchical
approach as proposed in the paper but by an {\em ad hoc} choice of $\pi
(\bm{\mu}|a)$.

%s5 ###
\section{Final remarks}

This paper contains many interesting ideas and examples. However, it
offers more of a brainstorming than a systematic treatment and a
general solution to the problem. It is somewhat disappointing that the
methods proposed in the paper bear little resemblance with the original
reference prior approach, where the problem is clearly stated, the
criterion used is sensible, and one can typically obtain \emph{unique}
and reasonable solutions. The approaches proposed here are still far
from becoming operational algorithms since they require a number of
arbitrary inputs. Hopefully, at least one of these methods will evolve
into an \emph{overall objective procedure} to find overall objective
priors. We believe the reference distance method to be the most
promising in this regard.

\bibliographystyle{ba}

\end{document}